\documentclass[reqno,12pt]{amsart}
\usepackage{url}
\usepackage{amssymb,amsthm,amsmath,amsfonts,amstext,mathrsfs}
\usepackage{epsfig}
\usepackage{latexsym}
\usepackage{fancyhdr}
\usepackage[unicode]{hyperref}

\hypersetup{colorlinks=true,
linkcolor=blue,
citecolor=blue,
urlcolor=blue,
filecolor=blue,
bookmarksnumbered=true,
pdfstartview=FitH,
pdfhighlight=/N
}

\newtheorem{thm}{Theorem}

\newtheorem{lem}{Lemma}

\newtheorem{defin}{Definition}

\newcommand{\p}{\partial}

\newcommand{\m}{\mathbf}
\newcommand{\mf}{\mathfrak}

\title[Giedrius Alkauskas]
{Beltrami vector fields with polyhedral symmetries} 
\author[Giedrius Alkauskas]{Giedrius Alkauskas}
\address{Vilnius University, Department of Mathematics and Informatics, Naugarduko 24, LT-03225 Vilnius, Lithuania}
\email{giedrius.alkauskas@mif.vu.lt}

\begin{document}
\begin{abstract} A $3$-dimensional vector field $B$ is said to be \emph{Beltrami vector field} (force free-magnetic vector field in physics), if $B\times(\nabla\times B)=0$. Motivated by our investigations on projective an polynomial superflows, and as an important side result, in the first paper on this topic we constructed two unique Beltrami vector fields $\mf{I}$ and $\mf{Y}$, such that $\nabla\times\mf{I}=\mf{I}$, $\nabla\times\mf{Y}=\mf{Y}$, and that both have orientation-preserving icosahedral symmetry (group of order $60$).\\
\indent In the current paper we extend these results to the tetrahedral and octahedral cases, and (together with an icosahedral case) we  calculate all simplest Beltrami fields with polyhedral symmetries arising from solutions to the Helmholtz equation of any order (the first aforementioned paper being an order $1$ approach).\\
\indent The notion of Beltrami vector field, slightly relaxed, generalizes to any dimension. In this paper we also present $2$-dimensional vector fields which have a dihedral symmetry $\mathbb{D}_{2d+1}$ of order $4d+2$. A much more detailed analysis is carried out in case $d=1$. One of these fields is particularly exceptional since it is the only case in our investigations which arises from the order $0$ approach to the Helmholtz equation, thus relating this flow to the $ABC$ flow.        
\end{abstract}

\pagestyle{fancy}
\fancyhead{}
\fancyhead[LE]{{\sc Beltrami vector fields}}
\fancyhead[RO]{{\sc G. Alkauskas}}
\fancyhead[CE,CO]{\thepage}
\fancyfoot{}

\date{\today}
\thanks{The research of the author was supported by the Research Council of Lithuania grant No. MIP-072/2015}

\maketitle

\section{Introduction}
\subsection{Groups}Let us define \small
\begin{eqnarray}
\alpha\mapsto\begin{pmatrix}
1 & 0 & 0\\
0 & -1 & 0\\
0 & 0 & -1 
\end{pmatrix},\quad
\beta\mapsto\begin{pmatrix}
0 & 1 & 0\\
1 & 0 & 0\\
0 & 0 & 1 
\end{pmatrix},
\quad
\gamma\mapsto\begin{pmatrix}
0 & 1 & 0\\
0 & 0 & 1\\
1 & 0 & 0 
\end{pmatrix},
\quad
\delta\mapsto\begin{pmatrix}
0 & 1 & 0\\
1 & 0 & 0\\
0 & 0 & -1 
\end{pmatrix}.
\label{g-24}
\end{eqnarray}\normalsize
Thus,
\begin{eqnarray*}
\langle\alpha,\gamma\rangle&=&\mathbb{T},\quad|\mathbb{T}|=12,\\
\langle\alpha,\beta,\gamma\rangle&=&\widehat{\mathbb{T}},\quad |\widehat{\mathbb{T}}|=24,\\
\langle\alpha,\gamma,\delta\rangle&=&\mathbb{O},\quad |\mathbb{O}|=24.
\end{eqnarray*}
The first is the tetrahedral group, the second - the full tetrahedral group, and the last one - the octahedral group. The Klein $4$-group $\mathbb{K}$, given by matrices 
\begin{eqnarray*}
\{\mathrm{diag}(1,-1,-1),\mathrm{diag}(-1,1,-1),\mathrm{diag}(-1,-1,1),I\}
\end{eqnarray*}
is a subgroup of $\mathbb{T}$.\\

Next, let
\begin{eqnarray*}
\eta\mapsto\begin{pmatrix}
\frac{1}{2} & -\frac{\phi}{2} & \frac{1}{2\phi}\\
\frac{\phi}{2} & \frac{1}{2\phi} & -\frac{1}{2}\\
\frac{1}{2\phi} & \frac{1}{2} & \frac{\phi}{2} 
\end{pmatrix},\quad \langle\mathbb{T},\eta\rangle=\mathbb{I},\quad |\mathbb{I}|=60.
\end{eqnarray*}
$\mathbb{I}$ is the icosehedral group.
\subsection{Beltrami vector fields}We quickly remind the method we used to construct Baltrami vector fields.\\

One of the main identities of the vector calculus claims that, for smooth vector field $\m{B}$, one has
\begin{eqnarray} 
\nabla\times(\nabla\times\m{B})=\nabla(\nabla\cdot\m{B})-\nabla^{2}\m{B}.
\label{vec-bas}
\end{eqnarray} 
Suppose now, a vector field $\m{B}$ satisfies the \emph{vector Helmholtz equation}
\begin{eqnarray*}
\nabla^{2}\m{B}=-\m{B},
\end{eqnarray*}
and also $\nabla\cdot\m{B}=0$. Then the identity (\ref{vec-bas}) gives
\begin{eqnarray*}
\nabla\times(\nabla\times\m{B})=\m{B}.
\end{eqnarray*}
Therefore, if $\nabla\times\m{B}=\m{C}$, then
\begin{eqnarray*}
\nabla\times(\m{B}+\m{C})=\m{B}+\m{C},
\end{eqnarray*}
and therefore $\m{B}+\m{C}$ is a Beltrami field.
\subsection{Order $n$ approach via a Helmholtz equation}
Let, as usual, 
\begin{eqnarray*}
P_{n}(x,y)=\Re\big{(}(x+iy)^{n}\big{)},\quad Q_{n}(x,y)=\Im\big{(}(x+iy)^{n}\big{)}
\end{eqnarray*}
be the standard harmonic polynomials of order $n$. These satisfy
\begin{eqnarray*}
\frac{\p P_{n}}{\p x}=nP_{n-1},&\quad& \frac{\p P_{n}}{\p y}=-nQ_{n-1},\\
\frac{\p Q_{n}}{\p x}=nQ_{n-1},&\quad& \frac{\p Q_{n}}{\p y}=nP_{n-1}.
\end{eqnarray*}
We will construct Beltrami vector fields from the order $n$ solutions to the Helmholtz equation which are given by the following lemma.
\begin{lem}
\label{lemu}
Let $\m{a}$, $\m{b}$ and $\m{c}$ be three $3$-dimensional vectors-rows, and $\m{x}=(x,y,z)^{T}$. The function $F(\m{x})=dP_{n}(\m{a}\m{x},\m{b}\m{x})\sin(\m{b}\m{x})$, $d\in\mathbb{R}$, is a solution to the Helmholtz equation if  $\m{a},\m{b}$ and $\m{c}$ are orthonormal. The same holds for the $\cos$ function, or the polynomial $Q_{n}$ instead of $P_{n}$.
\end{lem} 
\subsection{Lambent flows}
We proceed with a definition.
\begin{defin}
Let $\Gamma\hookrightarrow\mathrm{GL}(3,\mathbb{R})$ be an exact irreducible representation of a finite group ($\mathbb{T}$, $\mathbb{O}$ or $\mathbb{I}$). Suppose, a vector field $F$ satisfies these properties:
\begin{itemize}
\item[i)] $\nabla\times F=F$;
\item[ii)] All coordinates of $F$ are entire functions;
\item[iii)] if $\Gamma$ is the group in consideration, and if $F$ is treated as a map $\mathbb{R}^{3}\mapsto\mathbb{R}^{3}$, then $\epsilon^{-1}\circ F\circ\epsilon=F$ for any $\epsilon\in\Gamma$.
\end{itemize}
Then such a vector field is called \emph{lambent vector field with a symmetry} $\Gamma$.
\end{defin}
In dimension other that $3$ we do not have same dimensional analogue of the curl operator, but, as we have seen, it is natural to define lambent flow in any dimension as follows.  

\begin{defin}
Let $\Gamma\hookrightarrow\mathrm{GL}(n,\mathbb{R})$, $n\neq 3$, be an exact irreducible representation of a finite group. Suppose, a vector field $F$ satisfies these properties:
\begin{itemize}
\item[i)] $\nabla^{2}F=-F$;
\item[ii)]$\nabla\cdot F=0$;
\item[iii)] All coordinates of $F$ are entire functions and if expanded as Taylor series, they have only even compound degrees;
\item[iv)] if $\Gamma$ is the group in consideration, and if $F$ is treated as a map $\mathbb{R}^{n}\mapsto\mathbb{R}^{n}$, then $\epsilon^{-1}\circ F\circ\epsilon=F$ for any $\epsilon\in\Gamma$.
\end{itemize}
Then such a vector field is called \emph{lambent vector field with a symmetry} $\Gamma$.
\label{defin2}
\end{defin}
\section{The tetrahedral and octahedral cases}
In the next two section we will concentrate on a $3$-dimensional case.
\subsection{$n$ is odd.}Let $n\in\mathbb{N}$ be an odd positive integer. We will construct a vector field $\mathfrak{V}$ whose first coordinate $\mathfrak{G}$ is made out of linear combinations of the expressions of the form
\begin{eqnarray*}
T_{n}(x,y)\sin z,
\end{eqnarray*}
where $T$ is a harmonic polynomial, such that
\begin{itemize}
\item[i)]$\mf{G}$ is of even compound degree;
\item[ii)]$\mf{V}$ has a tetrahedral symmetry $\mathbb{T}$;
\item[iii)]$\mf{G}$ satisfies the Helmholtz equation.
\end{itemize}
For this we must put
\begin{eqnarray}
\mf{V}=\Big{(}\mf{G}(x,y,z),\mf{G}(y,z,x),\mf{G}(z,x,y)\Big{)}.
\label{v-vec}
\end{eqnarray}
We furter require
\begin{itemize}
\item[iv)]$\nabla\cdot\mf{V}=0$.
\end{itemize}
next, since $\mathbb{K}<\mathbb{T}$, $\mf{G}$ should be of even degree in $x$ and of odd in each of $y,z$, where, as usual, $\cos$ and $\sin$ are counted as even and odd, respectively. Thus, let
\begin{eqnarray*}
\mf{G}_{n;a,b}=aQ_{n}(x,y)\sin z+bQ_{n}(x,z)\sin y,\quad a,b\in\mathbb{R}.
\end{eqnarray*}
The vector field $\mf{V}$ now has a tetrahedral symmetry, all its components have an even compound degree, and it satisfies the vector Helmholtz equation. We are left to calculate its divergence.\\ 

Indeed,
\begin{eqnarray*}
\frac{\p}{\p x}\mf{G}_{a,b}(x,y,z)=anQ_{n-1}(x,y)\sin z+bnQ_{n-1}(x,z)\sin y.
\end{eqnarray*}
So, by a direct calculation,
\begin{eqnarray}
\mathrm{div}\,\mf{V}_{n;a,b}=n\big{(}a+(-1)^{\frac{n+1}{2}}b\big{)}\Big{(}
Q_{n-1}(x,y)\sin z+Q_{n-1}(y,z)\sin x+Q_{n-1}(z,x)\sin y\Big{)}.
\label{div-1}
\end{eqnarray}
Here we used an identity that, for $m$ even,
\begin{eqnarray*}
Q_{m}(y,x)=(-1)^{\frac{m+2}{2}}Q_{m}(x,y).
\end{eqnarray*}
\subsection{$n$ is even}Now, if $n$ is even (to avoid confusion, we thus take the first coordinate of the vector field to be $\mf{H}$, and vector field itself as $\mf{W}$), all vector fields of the form (\ref{v-vec}), which have a tetrahedral symmetry, and satisfy the Helmholtz equation, are of even degree in $x$, of odd in each of $y,z$, are given by
\begin{eqnarray*}
\mf{H}_{n;c}=cQ_{n}(y,z)\cos x.
\end{eqnarray*}
By a direct calculation,
\begin{eqnarray}
\mathrm{div}\,\mf{W}_{n;c}=-c\Big{(}Q_{n}(x,y)\sin z+Q_{n}(y,z)\sin x+Q_{n}(z,x)\sin z\Big{)}.
\label{div-2}
\end{eqnarray}
\subsection{Solenoidality}Inspecting (\ref{div-1}) and (\ref{div-2}), we see that if 
\begin{eqnarray*}
\sum\limits_{n\text{ odd}}\mf{V}_{n;a,b}+\sum\limits_{n\text{ even}}\mf{W}_{n;c}
\end{eqnarray*}
is a solenoidal vector field, then it as a linear combinations of vector fields, whose first coordinate is
\begin{eqnarray}
\mf{G}_{2m+1;a,b}+\mf{H}_{2m;c},\text{ where }(2m+1)a+(-1)^{m+1}(2m+1)b-c=0.
\label{m-cond}
\end{eqnarray}
Thus, we have two independent solutions. We get a particularly elegant pair of solutions by requiring that $\mf{G}_{2m+1;a,b}+\mf{H}_{2m;c}$ is symmetric with respect to $y,z$, and so the vector field has a full tetrahedral symmetry; or it is anti-symmetric, and so the vector field has an octahedral symmetry.
\subsection{The tetrahedral symmetry}Thus, in (\ref{m-cond}) we have two cases. First, consider $m=2\ell$, $\ell\in\mathbb{N}_{0}$. If $\mf{G}_{4\ell+1;a,b}$ is symmetric in $y,z$, then $a=b$, and the condition (\ref{m-cond}) tells that $c=0$. Therefore, we have a vector field, whose first coordinate is
\begin{eqnarray*}
\mf{G}=Q_{4\ell+1}(x,y)\sin z+Q_{4\ell+1}(x,z)\sin y.
\end{eqnarray*}
Consider $m=2\ell+1$, $\ell\in\mathbb{N}_{0}$. Still we have $a=b$, but then (\ref{m-cond}) gives $c=(8\ell+6)a$. This gives a vector field
\begin{eqnarray*}
\mf{G}=Q_{4\ell+3}(x,y)\sin z+Q_{4\ell+3}(x,z)\sin y+(8\ell+6)Q_{4\ell+2}(y,z)\cos x.
\end{eqnarray*}
Finally, if $\mf{V}$ is given by (\ref{v-vec}), then, as usually, $\mf{V}+\nabla\times\mf{V}$ is a Beltrami vector field with a tetrahedral symmetry. By calculating, we have the following result.
\begin{thm}Let $\ell\in\mathbb{N}_{0}$. Let us define
\begin{eqnarray*}
\mf{G}&=&Q_{4\ell+1}(x,y)\sin z+Q_{4\ell+1}(x,z)\sin y+Q_{4\ell+1}(z,x)\cos y-Q_{4\ell+1}(y,x)\cos z,
\end{eqnarray*}
and also
\begin{eqnarray*}
\mf{G}&=&Q_{4\ell+3}(x,y)\sin z+Q_{4\ell+3}(x,z)\sin y+(8\ell+6)Q_{4\ell+2}(y,z)\cos x\\
&+&Q_{4\ell+3}(z,x)\cos y-Q_{4\ell+3}(y,x)\cos z\\
&-&(8\ell+6)P_{4\ell+2}(y,z)\sin x\\
&+&(8\ell+6)(4\ell+2)P_{4\ell+1}(x,y)\cos z-(8\ell+6)(4\ell+2)P_{4\ell+1}(x,z)\cos y.
\end{eqnarray*}
Then the vector field, defined by (\ref{v-vec}), has a tetrahedral symmetry $\mathbb{T}$, and satisfies $\nabla\times\mf{V}=\mf{V}$.
\label{thm1}
\end{thm}
The first formula, when $\ell=0$ and $\ell=1$ gives, in particular,
\begin{eqnarray*}
\mf{G}&=&y\sin z+z\sin y+x\cos y-x\cos z,\\
\mf{G}&=&(5x^4y-10x^2y^3+y^5)\sin z+(5x^4z-10x^2z^3+z^5)\sin y\\
&+&(5z^4x-10z^2x^3+x^5)\cos y-(5y^4x-10y^2x^3+x^5)\cos z.
\end{eqnarray*}
The second formula for $\ell=0$ gives
\begin{eqnarray*}
\mf{G}&=&(3x^2y-y^3)\sin z+(3x^2z-z^3)\sin y+12yz\cos x\\
&+&(3z^2x-x^3)\cos y-(3y^2x-x^3)\cos z\\
&+&6(z^2-y^2)\sin x+12x\cos z-12x\cos y.
\end{eqnarray*}
In the last example we easily calculate that, indeed, $\nabla\times\mf{V}=\mf{V}$.
\subsection{The octahedral symmetry}Consider again $m=2\ell+1$. If $\mf{G}_{4\ell+3;a,b}$ is antisymmetric in $a,b$, then $a=-b$, and then $c=0$. We thus again have the first coordinate
\begin{eqnarray*}
\mf{G}=Q_{4\ell+3}(x,y)\sin z-Q_{4\ell+3}(x,z)\sin y.
\end{eqnarray*}
Consider now $m=2\ell$. Still we have $b=-a$, but then $c=(8\ell+2)a$. This gives
\begin{eqnarray*}
\mf{G}=Q_{4\ell+1}(x,y)\sin z-Q_{4\ell+1}(x,z)\sin y+(8\ell+2)Q_{4\ell}(y,z)\cos x.
\end{eqnarray*}

 Again, calculating $\mf{V}+\nabla\times\mf{V}$, we have
\begin{thm}Let $\ell\in\mathbb{N}_{0}$. Let us define
\begin{eqnarray*}
\mf{G}&=&Q_{4\ell+3}(x,y)\sin z-Q_{4\ell+3}(x,z)\sin y+Q_{4\ell+3}(z,x)\cos y+Q_{4\ell+3}
(y,x)\cos z,
\end{eqnarray*}
and also
\begin{eqnarray*}
\mf{G}&=&Q_{4\ell+1}(x,y)\sin z-Q_{4\ell+1}(x,z)\sin y+(8\ell+2)Q_{4\ell}(y,z)\cos x\\
&+&Q_{4\ell+1}(z,x)\cos y+Q_{4\ell+1}(y,x)\cos z\\
&-&(8\ell+2)P_{4\ell}(y,z)\sin x\\
&+&4\ell(8\ell+2)P_{4\ell-1}(x,y)\cos z+4\ell(8\ell+2)P_{4\ell-1}(x,z)\cos y.
\end{eqnarray*}
Then the vector field, defined by (\ref{v-vec}), has an octahedral symmetry $\mathbb{O}$ (also, a tetrahedral one), and satisfies $\nabla\times\mf{V}=\mf{V}$.
\end{thm}
For $\ell=0$, the first formula gives
\begin{eqnarray*}
\mf{G}=(3x^2y-y^3)\sin z-(3x^2z-z^3)\sin y+(3z^2x-x^3)\cos y+(3y^2x-x^3)\cos z,
\end{eqnarray*}
and the second formula gives $(Q_{0}=P_{-1}=Q_{-1}=0)$
\begin{eqnarray*}
\mf{G}=y\sin z-z\sin y+x\cos y-2\sin x+x\cos z.
\end{eqnarray*}
\section{Icosahedral symmetry}
\subsection{Induced vector fields}Now, $\mathbb{T}$ is not a normal subgroup of $\mathbb{I}$, but rather we have the identities
\begin{eqnarray*}
\eta\cdot\alpha\cdot\eta^2=\hat{\alpha}\in\mathbb{T},\quad \eta\cdot\gamma\cdot\eta^3=\hat{\gamma}\in\mathbb{T}.
\end{eqnarray*}
Let $\mf{V}$ be any lambent vector field for the group $\mathbb{T}$. Let us consider the induced vector field, given by
\begin{eqnarray*}
\mf{V}^{\mathbb{I}}_{\eta}=\sum\limits_{j=0}^{4}\eta^{-j}\circ\mf{V}\circ\eta^{j}.
\end{eqnarray*}
First, each summand satisfies the vector Helmholtz equation, minding the invariance of the Helmholtz equation under the orthogonal transformation; therefore, so does the whole sum.  Second,
\begin{eqnarray*}
\nabla\times(\eta^{-j}\circ\mf{V}\circ\eta^{j})=\eta^{-j}\circ\mf{V}\circ\eta^{j}.
\end{eqnarray*}
Third, obviously, $\mf{V}^{\mathbb{I}}_{\eta}$ is invariant under conjugation with $\eta$. Finally, it is invariant under conjugation with $\alpha$ and $\gamma$, since 
\begin{eqnarray*}
\alpha^{-1}\circ\eta^{-1}\circ\mf{V}\circ\eta\circ\alpha=\eta^{-3}\circ\hat{\alpha}^{-1}\circ\mf{V}
\circ\hat{\alpha}\circ\eta^{3}&=&\eta^{-3}\circ\mf{V}\circ\eta^{3},\\
\gamma^{-1}\circ\eta^{-1}\circ\mf{V}\circ\eta\circ\gamma=\eta^{-2}\circ\hat{\gamma}^{-1}\circ\mf{V}
\circ\hat{\gamma}\circ\eta^{2}&=&\eta^{-2}\circ\mf{V}\circ\eta^{2}.
\end{eqnarray*}
More generally, we find that
\begin{eqnarray*}
\alpha^{-1}\circ\eta^{-j}\circ\mf{V}\circ\eta^{j}\circ\alpha&=&\eta^{-\sigma(j)}\circ\mf{V}\circ\eta^{\sigma(j)},
\quad \sigma(0,1,2,3,4)=(0,3,4,1,2),\\
\gamma^{-1}\circ\eta^{-j}\circ\mf{V}\circ\eta^{j}\circ\gamma&=&\eta^{-\phi(j)}\circ\mf{V}\circ\eta^{\phi(j)},\quad
\phi(0,1,2,3,4)=(0,2,4,3,1).
\end{eqnarray*}
Of course, we can start just from the even part of the vector field. For example, let us take the even part of $\mf{G}_{1}$ (using the notation of Theorem \ref{thm1}), which is just 
\begin{eqnarray*}
\mf{G}=y\sin z+z\sin y.
\end{eqnarray*}
By a direct calculation, computing the induced vector field, we obtain exactly the vector field given in \cite{alkauskas-ico}. We can reproduce it here:
\scriptsize
\begin{eqnarray*}
\mf{V}^{\eta}_{I}&=&2x\sin\Big{(}\frac{x}{2}\Big{)}\sin\Big{(}\frac{\phi y}{2}\Big{)}\sin\Big{(}\frac{z}{2\phi}\Big{)}
-2\phi x\sin\Big{(}\frac{x}{2\phi}\Big{)}\sin\Big{(}\frac{y}{2}\Big{)}\sin\Big{(}\frac{\phi z}{2}\Big{)}+2\phi^{-1} x\sin\Big{(}\frac{\phi x}{2}\Big{)}\sin\Big{(}\frac{y}{2\phi}\Big{)}\sin\Big{(}\frac{z}{2}\Big{)}\\
&+&y\sin z+2y\cos\Big{(}\frac{x}{2}\Big{)}\cos\Big{(}\frac{\phi y}{2}\Big{)}\sin\Big{(}\frac{z}{2\phi}\Big{)}
-2y\cos\Big{(}\frac{x}{2\phi}\Big{)}\cos\Big{(}\frac{y}{2}\Big{)}\sin\Big{(}\frac{\phi z}{2}\Big{)}\\
&+&z\sin y-2z\cos\Big{(}\frac{x}{2}\Big{)}\sin\Big{(}\frac{\phi y}{2}\Big{)}\cos\Big{(}\frac{z}{2\phi}\Big{)}+2z\cos\Big{(}\frac{\phi x}{2}\Big{)}\sin\Big{(}\frac{y}{2\phi}\Big{)}\cos\Big{(}\frac{z}{2}\Big{)}.
\end{eqnarray*}
\normalsize
If we alternatively put
\begin{eqnarray*}
\mf{G}^{0}=\phi y\sin z-\phi^{-1}z\sin y,
\end{eqnarray*}
we obtain a vector field in Theorem 2 in \cite{alkauskas-ico}. Of course, any vector field $a y\sin z+bz\sin y$ will do, but $\mf{G}$ and $\mf{G}^{0}$ have an additional symmetry with respect to a non-trivial automorphism $\tau$ of $\mathbb{Q}[\,\sqrt{5}\,]$. Indeed,
\begin{eqnarray*}
\tau\mf{G}(z,y)=\mf{G}(y,z),\quad \tau\mf{G}^{0}(z,y)=\mf{G}^{0}(y,z).
\end{eqnarray*}
We remind that $\tau$ acts on entire functions with rational Taylor coefficients the same way it acts on $\mathbb{Q}$; that is, leaves the coefficients (and the function) intact.

\section{The case $\mathbb{R}^{2}$}
\label{sec-2dim}
Now let us turn our attention to a $2$-dimensional case, based on Definition  \ref{defin2}. Let, as before, $P_{n}$ and $Q_{n}$ stand for the harmonic polynomials of order $n$. As we know from \cite{alkauskas-super1}, for $d\in\mathbb{N}$, a vector field
\begin{eqnarray*}
\Big{(}P_{2d}+(-1)^{d}Q_{2d},(-1)^{d}P_{2d}-Q_{2d}\Big{)}
\end{eqnarray*} 
has a dihedral symmetry $\mathbb{D}_{2d+1}$ of order $4d+2$, and is solenoidal. Thus, potentially it gives rise to a lambent flow. We will explore the case $d=1$, and only make the first few steps in a research, since the vector fields we obtain already contain very rich dynamic and analytic structure.\\

Thus, let (we multiply the above vector field by a factor $-1$ just for convenience)
\begin{eqnarray*}
(\varpi,\varrho)=(2xy-x^2+y^2,2xy+x^2-y^2),
\end{eqnarray*}
which is solenoidal, and has a dihedral group $\mathbb{D}_{3}$, generated by matrices 
\begin{eqnarray*}
\alpha=\frac{1}{2}\begin{pmatrix} -1 & -\sqrt{3}\\ 
\sqrt{3} & -1\\ \end{pmatrix},\quad
\beta=\begin{pmatrix} 0 & 1\\ 1 & 0\\ \end{pmatrix},
\end{eqnarray*} 
as its symmetry group. This gives rise to the dihedral superflow. The vector field, linearly conjugate (over $\mathrm{GL}_{2}(\mathbb{R})$, not over $O(2)$) to $(\varpi,\varrho)$ is given by $(x^2-2xy,y^2-2xy)$. This vector field is indeed a very fascinating one, and it was explored in detail in \cite{alkauskas-un}. For example, this vector field can be explicitly integrated in terms Dixonian elliptic functions  \cite{flajolet-c}. In particular, this flow turns out to be \emph{unramified}. The property of the flow to be unramified is indeed a profound one. For example, to prove that a flow with a vector field $(x^2-xy,y^2-2xy)$ is unramified, requires to employ arithmetic of the number field $\mathbb{Q}(\,\sqrt{3}\,)$ \cite{alkauskas-ab}.  The very theory of unramified flows is still in its initial stages of development.\\

 Returning to Beltrami flows in dimension $2$ and to the group $\mathbb{D}_{3}$ (analogously as in the icosahedral case in \cite{alkauskas-ico}), let us define $4$ linear forms and $4$ vectors
\begin{eqnarray*}
\begin{tabular}{l l}
$\ell_{x,0}(x,y)=\frac{1}{2}x+\frac{\sqrt{3}}{2}y$,& $\m{j}_{x,0}=\Big{(}\frac{1}{2},\frac{\sqrt{3}}{2}\Big{)}$,\\
$\ell_{x,1}(x,y)=-\frac{1}{2}x+\frac{\sqrt{3}}{2}y$,& $\m{j}_{x,1}=\Big{(}-\frac{1}{2},\frac{\sqrt{3}}{2}\Big{)}$,\\
$\ell_{y,0}(x,y)=\frac{\sqrt{3}}{2}x+\frac{1}{2}y$,& $\m{j}_{y,0}=\Big{(}\frac{\sqrt{3}}{2},\frac{1}{2}\Big{)}$,\\
$\ell_{y,1}(x,y)=\frac{\sqrt{3}}{2}x-\frac{1}{2}y$,& $\m{j}_{y,1}=\Big{(}\frac{\sqrt{3}}{2},-\frac{1}{2}\Big{)}$.
\end{tabular}
\end{eqnarray*}
We have a vector field $(\mf{G}(x,y),\mf{G}(y,x))$, where
\begin{eqnarray}
\mf{G}&=&a\cos x+b\cos y+c\cos\ell_{x,0}+d\cos\ell_{x,1}
+e\cos\ell_{y,0}+f\cos\ell_{y,1}\nonumber\\
&+&gy\sin x+hx\sin y\\
&+&i\ell_{y,1}\sin\ell_{x,0}+j\ell_{y,0}\sin\ell_{x,1}\nonumber\\
&+&k\ell_{x,1}\sin\ell_{y,0}+l\ell_{x,0}\sin\ell_{y,1}.\nonumber
\end{eqnarray}
We want the vector field $(\mf{G}(x,y),\mf{G}(y,x))$ to have a dihedral symmetry, to satisfy the Helmholtz equation, and to be  solenoidal. We also require that such a field has $(\varpi,\varrho)$ as its degree $2$ term in its Taylor series expansion. Thus, we have $12$ free parameters. Now, let us go through the same procedure as before with a help of MAPLE, what was done in \cite{alkauskas-ico} for the icosahedral case.\\

\indent\textbf{i)} The requirement that a Taylor series for $\mf{G}$ starts at $\varpi$ puts three linear conditions, the dependent parameters being $f,l,k$.\\
\indent\textbf{ii)} Solenoidality of $(\mf{G}(x,y),\mf{G}(y,x))$ puts $5$ further requirements, leaving $a,b,c,d$ as free coefficients.\\
\indent\textbf{iii)} Finally, the invariance under $\alpha$ requires that
\begin{eqnarray*}
-\frac{1}{2}\mf{G}(-\ell_{x,0},\ell_{y,1})+
\frac{\sqrt{3}}{2}\mf{G}(\ell_{y,1},-\ell_{x,0})=\mf{G}(x,y).
\end{eqnarray*} 
This gives two more conditions, leaving $a$ as a free parameter, thus giving a $1$-parameter family which solves our problem. Two independent choices choices of $a$ that give a particularly symmetric solututions turn out to be choose $a=0$ and $a=\frac{4}{3}$.\\

In the first case, the whole collection turn out to be
\begin{eqnarray*}
(a,b,c,d,e,f,g,h,i,j,k,l)=\big{(}0,-\frac{8}{3},-\frac{4\sqrt{3}}{3},
\frac{4\sqrt{3}}{3},\frac{4}{3},\frac{4}{3},0,0,0,0,0,0\big{)}.
\end{eqnarray*}
Thus, we have an order $0$ solution. For convenience, we multiply all coefficient by $\frac{3}{8}$ throughout. This gives the formulas for the vector field, given by Theorem \ref{teom3}.\\

%\begin{figure}
%\epsfig{file=5orbits,width=320pt,height=350pt,angle=-90}
%\caption{Orbits for the selected $5$ points $(6,0)$, $(6.5,0)$, $(7.0,0)$, $(7.5,0)$ and %$(11.0,0)$.}
%\label{figure-dih}
%\end{figure}

\begin{figure}
\epsfig{file=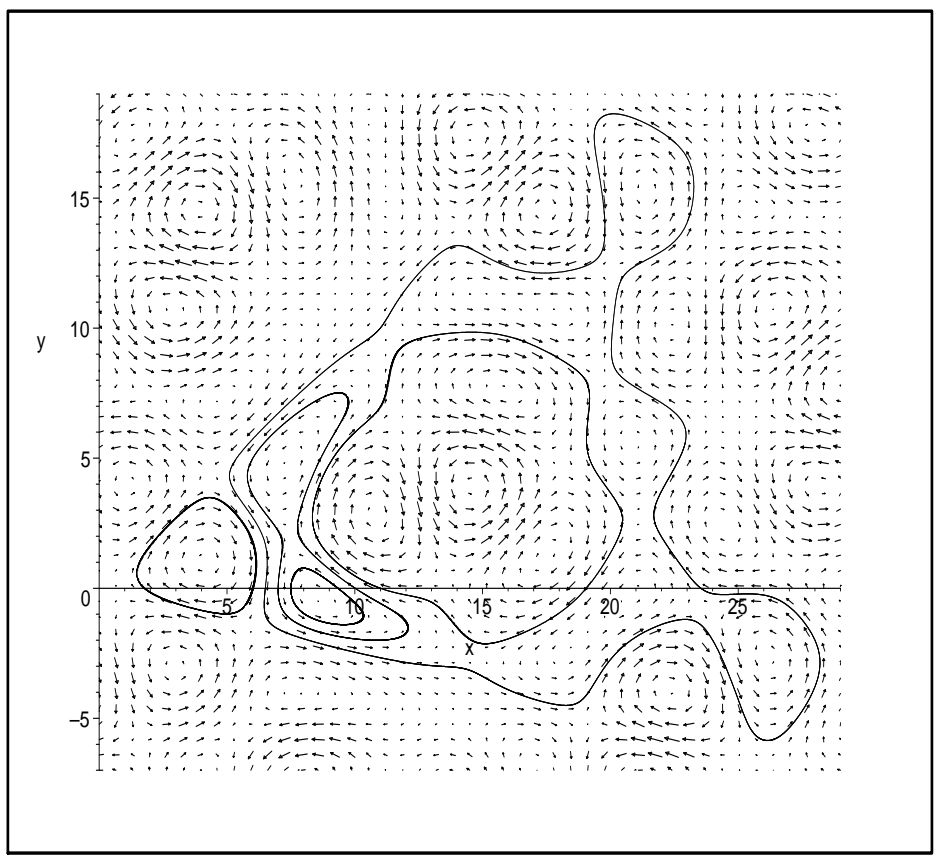,width=320pt,height=350pt,angle=-90}
\caption{Orbits for the selected $5$ points $(6,0)$, $(6.5,0)$, $(7.0,0)$, $(7.5,0)$ and $(11.0,0)$.}
\label{figure-dih}
\end{figure}

In Figure \ref{figure-dih}, we plot orbits for the selected $5$ points $(6,0)$, $(6.5,0)$, $(7.0,0)$, $(7.5,0)$ and $(11.0,0)$. The picture itself suggests the following remarks:
\begin{itemize}
\item[1)] Let the solutions of the transcendental equation $\varpi(x,x)=0$, $x\geq 0$ by given by $\tau_{i}$, $i\in\mathbb{N}$. The vector field vanishes for $(x,y)=(\tau_{i},\tau_{i})$. For other points on the line $x=y$, the flow $F(\m{x},t)$ is a flow on the line, satisfying $F(\m{x},-\infty)=(\xi,\xi)$, $F(\m{x},\infty)=(\chi,\chi)$, where $\xi$ and $\chi$ are two adjacent values of $\tau$; 
\item[2)] There exists a point on the line $y=0$, somewhere between $(6.0,0)$ and $(6.5,0)$ such that its orbits is an unbounded, or not-closed curve.
\end{itemize}  
The choice $a=\frac{4}{3}$, if we additionally multiply all coefficients by $\frac{3}{2}$ througout, gives the collection of parameters
\begin{eqnarray*}
(a,b,c,d,e,f,g,h,i,j,k,l)=\big{(}2,2,\sqrt{3}-1,-\sqrt{3}-1,-\sqrt{3}-1,\sqrt{3}-1,0,2,-\sqrt{3},\sqrt{3},1,1\big{)}.
\end{eqnarray*}
Thus we have the following.
\begin{thm}
\label{teom3}
Let us define the vector field $\mf{V}=(\mf{V}_{x},\mf{V}_{y})$, where
\begin{eqnarray*}
\mf{V}_{x}&=&-\cos y+\sqrt{3}\sin\Big{(}\frac{x}{2}\Big{)}\sin\Big{(}\frac{\sqrt{3}y}{2}\Big{)}+\cos\Big{(}\frac{\sqrt{3}x}{2}\Big{)}\cos\Big{(}\frac{y}{2}\Big{)},\\
\mf{V}_{y}&=&-\cos x+\sqrt{3}\sin\Big{(}\frac{y}{2}\Big{)}\sin\Big{(}\frac{\sqrt{3}x}{2}\Big{)}+\cos\Big{(}\frac{\sqrt{3}y}{2}\Big{)}\cos\Big{(}\frac{x}{2}\Big{)}.
\end{eqnarray*}
Also,
\begin{eqnarray*}
\mf{Q}_{x}&=&2\cos x+2\cos y+2x\sin y+2\sqrt{3}\sin\Big{(}\frac{\sqrt{3}x}{2}\Big{)}\sin\Big{(}\frac{y}{2}\Big{)}-x\cos\Big{(}\frac{\sqrt{3}x}{2}\Big{)}\sin\Big{(}\frac{y}{2}\Big{)}\\
&+&\sqrt{3}y\sin\Big{(}\frac{\sqrt{3}x}{2}\Big{)}\cos\Big{(}\frac{y}{2}\Big{)}-2\cos\Big{(}\frac{x}{2}\Big{)}\cos\Big{(}\frac{\sqrt{3}y}{2}\Big{)}+
\sqrt{3}y\cos\Big{(}\frac{x}{2}\Big{)}\sin\Big{(}\frac{\sqrt{3}y}{2}\Big{)}\\
&-&2\sqrt{3}\sin\Big{(}\frac{x}{2}\Big{)}\sin\Big{(}\frac{\sqrt{3}y}{2}\Big{)}-3x\sin\Big{(}\frac{x}{2}\Big{)}\cos\Big{(}\frac{\sqrt{3}y}{2}\Big{)}-2\cos\Big{(}\frac{\sqrt{3}x}{2}\Big{)}\cos\Big{(}\frac{y}{2}\Big{)},
\end{eqnarray*}
$\mf{Q}=(\mf{Q}_{x},\mf{Q}_{y})$, $\mf{Q}_{y}(x,y)=\mf{Q}_{x}(y,x)$. 
Then
\begin{itemize}
\item[1)] The vector field $\mf{V}$ and has a dihedral symmetry of order $6$: for any $\zeta\in\mathbb{D}_{3}$, one has $\zeta^{-1}\circ\mf{V}\circ\zeta=\mf{V}$;
\item[2)] as a Taylor series, $\mf{V}$ contains only terms with even compound degrees, and it satisfies the vector Helmholtz equation and solenoidality:
\begin{eqnarray*}
\nabla^{2}\mf{V}=-\mf{V},\quad \nabla\cdot\mf{V}=0;
\end{eqnarray*} 
both these properties hold for $\mf{Q}$, too;
\item[3)] the Taylor series for $\mf{V}$ and $\mf{Q}$ start form a degree $2$ vector field $\frac{3\m{M}}{8}$ and $\frac{3\m{M}}{2}$, resepctively, where $\m{M}=(\varrho,\varpi)$ is given 
\begin{eqnarray*}
(\varpi,\varrho)=(2xy-x^2+y^2,2xy+x^2-y^2).
\end{eqnarray*}
\end{itemize} 
\end{thm}

\end{document}